\documentclass[%
 aip,
rsi,%
 amsmath,amssymb,
]{revtex4-1}

\usepackage{graphicx}
\usepackage{dcolumn}
\usepackage{bm}

\usepackage{graphicx}
\usepackage{amsmath}
\usepackage{amssymb}
\usepackage{algorithm}
\usepackage[update,prepend]{epstopdf}
\usepackage{subfigure}
\usepackage{amsthm}
\usepackage{enumerate}
\usepackage{color}
\usepackage{yfonts}
\usepackage{siunitx}

\usepackage[version=4]{mhchem}

\usepackage[framemethod=tikz]{mdframed}

\def\fatr{\mathbf{r}}
\def\fatR{\mathbf{R}}

\def\kme{k_{\rm CK}}

\def\krmi{k_a}
\def\krme{\krmi^{\rm meso}}

\def\rrad{\sigma}
\def\diffc{D}

\def\microbind{\tau_{\rm{micro}}}
\def\mesobind{\tau_{\rm{meso}}}
\def\microdiff{\tau_{\rm{diff}}^{\rm{micro}}}
\def\mesodiff{\tau_{\rm{diff}}^{\rm{meso}}}

\def\microreact{\tau_{\rm{react}}^{\rm{micro}}}
\def\mesoreact{\tau_{\rm{react}}^{\rm{meso}}}
\def\hybridreact{\tau_{\rm{react}}^{\rm{hybrid}}}

\def\hstar{h^*}

\begin{document}

\title{Mesoscopic-microscopic spatial stochastic simulation with automatic system partitioning}

\author{Stefan Hellander}
 \affiliation{Department of Information Technology, Uppsala University, Box 337, SE-75105, Uppsala, Sweden.}
\author{Andreas Hellander}
\affiliation{Department of Information Technology, Uppsala University, Box 337, SE-75105, Uppsala, Sweden.}
\author{Linda Petzold}
\affiliation{Department of Computer Science, University of California,Santa Barbara, CA 93106-5070 Santa Barbara, USA.}

\date{\today}

\begin{abstract}
The reaction-diffusion master equation (RDME) is a model that allows for efficient on-lattice simulation of spatially resolved stochastic chemical kinetics. Compared to off-lattice hard-sphere simulations with Brownian Dynamics (BD) or Green's Function Reaction Dynamics (GFRD) the RDME can be orders of magnitude faster if the lattice spacing can be chosen coarse enough. However, strongly diffusion-controlled reactions mandate a very fine mesh resolution for acceptable accuracy. It is common that reactions in the same model differ in their degree of diffusion control and therefore require different degrees of mesh resolution. This renders mesoscopic simulation inefficient for systems with multiscale properties. Mesoscopic-microscopic hybrid methods address this problem by resolving the most challenging reactions with a microscale, off-lattice simulation. However, all methods to date require manual partitioning of a system, effectively limiting their usefulness as 'black-box' simulation codes. In this paper we propose a hybrid simulation algorithm with automatic system partitioning based on indirect \emph{a priori} error estimates. We demonstrate the accuracy and efficiency of the method on models of diffusion-controlled networks in 3D.
%
\end{abstract}

\pacs{Valid PACS appear here}
\maketitle

\section{\label{sec:intro}Introduction}

Cellular control systems are inherently spatial, as reactions in a network involve macromolecules confined to certain locations or sub-compartments. For example, MAPK pathways involve sensing a stimulus by receptors localized to the cell membrane, propagating the signal in the cytoplasm via a phosphorylation cascade that modifies transcription factors, which are eventually imported into the nucleus where they bind to promoter sites on DNA and affect the expression of downstream genes. On an even finer level, spatial correlations on very short length scales impact the dynamics of various biochemical systems \cite{Mahmutovic2012}. The low copy numbers of key molecules introduces stochasticity in gene regulatory networks (GRNs). This is an important factor to account for when studying the regulatory properties of GRNs. Examples where both spatial and stochastic effects are predicted to be important \cite{Mahmutovic2012} include spatial gene regulation of Hes1 \cite{Sturrock2:2013}, polarization in budding yeast \cite{Lawson:2013}, and the MinD-system in \emph{E. Coli} \cite{FaEl}.

As a consequence, spatio-temporal simulation of reaction-diffusion systems is an important tool to analyze GRNs. In particular, two modeling frameworks have attracted considerable attention in the systems biology community: the mesoscopic, discrete stochastic reaction-diffusion master equation (RDME) in which point-like molecules are tracked on a grid, and Brownian Dynamics (BD) in which hard-sphere particles are tracked in continuous space.  
Many capable software packages have been created to support such spatial modeling, including MCell \cite{mcell}, Smoldyn \cite{smoldyn}, E-Cell \cite{ecell}, MesoRD \cite{mesord}, VCell \cite{vcell}, STEPS \cite{steps}, NeuroRD \cite{neurord}, ReaDDy \cite{readdy}, URDME \cite{URDME_BMC}, PyURDME \cite{molns} and StochSS \cite{StochSS}, the latter which integrates spatial capabilities via PyURDME.

Mesoscopic simulators are efficient if a reasonably coarse mesh can be used. However, for some diffusion-limited systems, it is critical to capture short-range, short-timescale interactions between the molecules  \cite{TaTNWo10,FangeSRDME,Mahmutovic2012}. In this case a microscopic, particle-scale resolution is needed for accurate simulation. This raises the question of how well the mesoscopic model can capture the microscale dynamics as the mesh size approaches the molecule size. Unless  measures are taken, the RDME approximation quality degrades as the mesh size becomes increasingly fine \cite{doi:10.1137/070705039,HHP2,HP1}.

The point-particle mesoscopic model can be made to approximate a microscopic model by deriving scale-dependent reaction rates \cite{FangeSRDME,HHP,HHP2,ErCha09} down to a critical mesh size after which no one-neighbor stencil can provide increased accuracy \cite{HHP}. After this critical limit, it is possible to improve the approximation even further by considering a wider stencil \cite{HP1,FangeSRDME}. This approach results in a lattice method with a lattice spacing on the order of the size of the molecules. By considering the Doi model \cite{Doi1} rather than the Smoluchowski model on the microscopic scale, it is also possible to arrive at a non-local, convergent mesoscopic model by directly discretizing the microscopic model \cite{IsaacsonCRDME}. But even if these solutions allow for accurate simulation of diffusion-limited systems down to a mesh size close to individual particles, the simulation cost increases dramatically as the mesh becomes finer: the number of diffusion events (and so the simulation time) grows proportionally to $h^{-2}$ where $h$ is a measure of the mesh size.

A promising approach to efficiently simulate systems with multiscale properties are hybrid methods in which the reaction network is partitioned and parts of it are simulated on different modeling levels and with different solvers. Examples include mesoscopic-macroscopic (PDE) methods \cite{FERM2010343,tworegime1,SPILL2015429,Lo160485}, macroscopic-mesoscopic methods \cite{Yates20150141}, macroscopic-microscopic methods \cite{doi:10.1137/120882469}, and mesoscopic-microscopic methods, in which parts of a system are simulated with the RDME and parts of it are simulated with a particle-tracking algorithm \cite{hybrid1,hybrid2,doi:10.1093/bioinformatics/btv149}. If a good partitioning can be found, the cost savings of hybridization can be substantial by keeping the number of particles handled on the microscale small, yet maintaining a reasonable mesh size for the RDME simulator. However, one key problem with hybrid methods is that prior knowledge about the system is needed in order to partition it correctly. Also, the system dynamics may change over the course of the simulation or in different regions in the spatial domain, making the initial system partitioning invalid or suboptimal. These issues make hybrid solvers hard to use without expert knowledge, which limits their usefulness as black-box simulation tools. Another problem is that they are complex and challenging to implement, thus there is a lack of software capable of hybrid simulation.

In this paper we develop theory and a practical method that, given a user-supplied error tolerance, is capable of automatic selection of the appropriate modeling level for each species in a spatial stochastic model. We show that the hybrid method converges with decreasing time step, and demonstrate numerically that it accurately reproduces the fine-grained reaction dynamics of microscopic methods. Finally, we show that the hybrid method can simulate systems that are intractable with pure mesoscopic methods, without having to simulate the whole system microscopically. 

The rest of the paper is organized as follows. In Sect. \ref{sec:background} we briefly introduce the mesoscopic and microscopic models, and review how they are related. In Sect. \ref{sec:method} we describe the proposed hybrid algorithm, we show how to split a general system, and we develop the condition to be used for adaptive system partitioning.
In Sect. \ref{sec:implementation} we discuss the practical implementation of the method, and demonstrate the accuracy and efficiency of the developed method on two challenging test problems in Sect. \ref{sec:results}. Sect \ref{sec:conclusions} concludes the paper.

\section{Background}
\label{sec:background}
Two modeling frameworks that are popular for simulating reaction-diffusion kinetics in systems biology are the microscopic Smoluchowski model and the mesoscopic reaction-diffusion master equation (RDME). In the former, particles are modeled as  hard spheres and their positions are tracked continuously in space, whereas in the latter, particles are point-particles and their positions are tracked up to the resolution of a structured or unstructured grid approximating the domain.

\subsection{Mesoscopic model}
On the mesoscopic scale, we model molecules as point particles, and treat diffusion as jumps between adjacent voxels on a mesh. The state of the system is the discrete number $s_i$ of each chemical species $S_i$, $i=1\ldots M$ in the voxels of the grid, where the voxels are denoted by $\mathcal{V}_j,~j=1\ldots N$.

A diffusive jump is a linear event
\begin{align}
S_{1i} \xrightarrow{D_1} S_{1j},
\end{align}
taking a molecule of species $S_1$ from voxel $\mathcal{V}_i$ to an adjacent voxel $\mathcal{V}_j$, where $D_1$ is a rate constant that depends on the diffusion constant of species $S_1$ and the shape and size of the voxels \cite{efhl2009}.

Chemical reactions occur between molecules residing in the same voxel. For example, a bimolecular reaction between species $S_1$ and $S_2$ producing $S_3$ in voxel $\mathcal{V}_j$ can be written
\begin{align}
S_{1j} + S_{2j} \xrightarrow{k_a} S_{3j},
\end{align}
where $k_a$ is a mesoscopic reaction rate. The mathematical formalism is the continuous-time discrete-space Markov process. In this framework, the propensity for the reaction, $a(s_{1j},s_{2j}) = k_a s_{1j}s_{2j}/V_{j}$, is the probability of the reaction occurring in an infinitesimal interval $(t,t+\delta t)$, where $V_j$ is the volume of voxel $\mathcal{V}_j$. With this assumption, realizations of the process can be efficiently simulated using versions of the SSA optimized for reaction-diffusion systems such as the Next Subvolume Method (NSM) \cite{ElEh04}.

\subsection{Microscopic model}
\label{micro-background}

Consider two molecules $M_1$ and $M_2$, of species $S_1$ and $S_2$ respectively. The molecules can react irreversibly according to $\ce{S_1 + S_2 -> S_3}$. In the microscopic Smoluchowski model, molecules are modeled as hard spheres with a finite reaction radius, diffusing according to Brownian motion. We denote the molecules' reaction radii by $\sigma_1$ and $\sigma_2$, and their diffusion constants by $D_1$ and $D_2$. Their positions in $\mathbb{R}^3$ at time $t_0$ are denoted by $\fatr_{10}$ and $\fatr_{20}$. The probability that the molecules are at positions $\fatr_1$ and $\fatr_2$ at time $t$ is described by the probability density function (pdf) $p(\fatr_1,\fatr_2,t|\fatr_{10},\fatr_{20},t_0)$. Let $\fatr = \fatr_1-\fatr_2$, and $\fatR = \sqrt{\frac{D_2}{D_1}}\fatr_1+\sqrt{\frac{D_1}{D_2}}\fatr_2$. We can now rewrite the pdf as
\begin{align}
p(\fatr_1,\fatr_2,t|\fatr_{10},\fatr_{20},t_0) = p_{\fatR}(\fatR,t|\fatr_0,t_0)p_{\fatr}(\fatr,t|\fatr_0,t_0),
\end{align}
and it can be shown that the dynamics of the two molecules in $\mathbb{R}^3$ is governed by the following system of equations \cite{ZoWo5a}:
\begin{align}
\label{smolu-eq1}
\frac{\partial p_{\fatR}}{\partial t} = D\Delta_{\fatR} p_{\fatR}\\
\frac{\partial p_{\fatr}}{\partial t} = D\Delta_{\fatr} p_{\fatr},
\end{align}
where $D=D_1+D_2$. The initial and boundary conditions are given by
\begin{align}
4\pi\sigma^2D\frac{\partial p_{\fatr}}{\partial n}\bigg|_{\| \fatr\|=\rrad} &= \krmi p_{\fatr}(\| \fatr\| = \rrad, t|\fatr_0,t_0)\label{smolu-eq1-bicond1}\\
p_{\fatr}(\| \fatr\| \to \infty, t|\fatr_0,t_0) &= 0\label{smolu-eq1-bicond2}\\
p_{\fatr}(\fatr, t_0|\fatr_0,t_0) &= \delta(\fatr-\fatr_0).\label{smolu-eq1-bicond3}
\end{align}
This system of equations can be solved exactly \cite{ZoWo5a,CarJae}, or using an operator split approach \cite{SHeLo11}. We will call $\krmi$ the microscopic reaction rate. The probability that the two molecules react between times $t_0$ and $t$ is given by
\begin{align}
\label{micro_time}
p_{\fatr}(\ast,t|\fatr_0,t_0) = 1-\int_{t_0}^{t}p_{\fatr}(\fatr,t|\fatr_0,t_0)\,d\fatr.
\end{align}

The time $t_d$ until a molecule undergoes a unimolecular reaction event with reaction rate $k_d$ is given by sampling $t_d$ from an exponential distribution with mean $k_d$. If the dissociation event produces two molecules, then they are placed at a distance of $\rrad$ apart. We will outline in Sect. \ref{micro-implementation} how to use $p_{\fatr}(\ast,t|\fatr_0,t_0)$ and $p(\fatr_1,\fatr_2,t|\fatr_{10},\fatr_{20},t_0)$ to simulate a more complex system within a bounded domain.

\subsection{Mesoscopic parameters}
\label{mesoparams}
We now ask how the mesoscopic and microscopic models are related. Specifically we need to relate the mesoscopic parameters to the microscopic parameters. The diffusion jump rates on the mesoscopic scale are obtained by discretizing the diffusion equation. On a structured mesh this is straightforward; for details on how to do it on an unstructured mesh, see \cite{URDME_BMC}. Below we outline an approach to relating the mesoscopic reaction rates to the microscopic reaction rates.

It is well known that if the reaction volume $V$ in 3D is much larger than the molecules, i.e. $V\gg\rrad$, then the mesoscopic reaction rate, $\kme$, relates to the microscopic reaction rate as
\begin{align}
\label{ck-equation}
\kme = \frac{1}{V}\frac{4\pi\rrad D \krmi}{4\pi\rrad D + \krmi},
\end{align}  
where $D$ is the sum of the molecules' diffusion constants, and $\rrad$ is the sum of the reaction radii. This expression was first derived by Collins and Kimball \cite{CollinsKimball}, and later re-derived by Gillespie \cite{gillespierates}. It is easy to see that for a spatially discretized well-mixed system, or if the voxels are large enough ($h\gg\rrad$, where $h^3$ is the volume of a voxel), the mesoscopic reaction rate, $\krme$, will be given by
\begin{align}
\krme = \frac{1}{h^3}\frac{4\pi\rrad D\krmi}{4\pi\rrad D+\krmi}.
\end{align}
In models with multiscale properties we often need to resolve part of the system to very high accuracy, which requires a highly resolved mesh. This implies that the condition $h\gg\rrad$ might not be satisfied. We thus need to derive reaction rates for this case.

Following the analysis in \cite{HHP,HHP2}, we start by considering a single irreversible reaction
\begin{align}
\label{r-basic}
\ce{S_{1} + S_{2} ->[$\krmi$] S_3}
\end{align}
in a cube discretized by a Cartesian mesh. Additionally, we will assume that the $S_{1}$ molecule is fixed inside a voxel close to the center of the domain, while the $S_{2}$ molecule diffuses freely with diffusion rate $\diffc$.

To study the relationship between the RDME and Smoluchowski models for small $h$, we compare the expected time until the molecules first react on the microscopic scale to the time on the mesoscopic scale.

First, consider the microscopic scale. Let the $S_{2}$ molecule have an initial position sampled from a uniform distribution and denote the mean binding time, or the average time until the molecules react given that the $S_{1}$ molecule is uniformly distributed, by $\microbind$. Following the approach in \cite{HHP}, we split $\microbind$ into two parts:
\begin{align}
\microbind = \microdiff+\microreact,
\end{align}
where $\microdiff$ is the average time until the $S_{1}$ molecule is in contact with the $S_{2}$ molecule for the first time, and $\microreact$ is the average time until the molecules react given that they are in contact.

\noindent
We know that \cite{FBSE10,HHP2}:
\begin{align}
\microdiff \approx \begin{cases}
\frac{V}{4\pi\sigma D}, \quad (3D)\\
\frac{V\left\{ \log\left(\pi^{-1}\frac{V^{1/2}}{\sigma}\right)\right\}}{2\pi D}, \quad (2D)
\end{cases}
\label{eq:microdiff}
\end{align}
and that
\begin{align}
\microreact = \frac{V}{\krmi}\quad \text{(1D, 2D, 3D)}.
\label{eq:microreact}
\end{align}

We now consider the system \eqref{r-basic} on the mesoscopic scale. The $S_{2}$ molecule is fixed in a voxel close to the origin, so that it is far from the boundaries. The $S_{1}$ molecule is sampled uniformly on the mesh, and $\mesobind$ denotes the average time until the two molecules react for the first time. Let $\mesodiff$ be the average time until the molecules are in the same voxel for the first time, and let $\mesoreact$ denote the average time until they react given that the molecules start in the same voxel.

Again we split the average binding time into two parts
\begin{align}
\mesobind = \mesodiff+\mesoreact,
\end{align}
where (with $C_2=0.1951$ and $C_3=1.5164$) \cite{MoWe65,Montroll68,HHP,HHP2}
\begin{align}
\label{eq:mesodiff}
\mesodiff = \begin{cases}
\frac{C_3V}{6Dh}+O\left(N^{\frac{1}{2}}\right)\quad (3D)\\
\frac{V}{4\pi D}\log(N)+\frac{C_2 V}{4D}+O\left(N^{-1}\right)\quad (2D)
\end{cases}
\end{align}
and
\begin{align}
\label{eq:mesoreact}
\mesoreact = \frac{N}{\krme},
\end{align}
\noindent
We make the ansatz that the mean reaction times on the mesoscopic and microscopic scales are equal, to obtain
\begin{align}
\mesobind = \microbind \\
\iff \mesodiff+\mesoreact = \microbind \\
\iff \krme = \frac{N}{\microbind-\mesodiff},\label{eq:mesorate}
\end{align}
where the last equality follows from \eqref{eq:mesoreact}. Note that $\microbind$ and $\mesodiff$ are both known, so that we can compute $\krme$ using \eqref{eq:mesorate}. We showed in \cite{HHP2} that $\krme$ can be rewritten as
\begin{align}
\label{rate-eqs}
\krme = \frac{\krmi}{h^3}\left( 1+\frac{\krmi}{D}G(h,\rrad) \right)^{-1}
\end{align}
where $G$ in 3D is given by
\begin{align}
G(h,\rrad) = \frac{1}{4\pi\rrad}-\frac{C_3}{6h}.
\end{align}

From the analysis follows the existence of an (for accuracy) optimal mesh size. To see this, consider that:
\begin{align}
\mesodiff < \microdiff &\implies \mesoreact > \microreact\label{dyn1}\\
\mesodiff = \microdiff &\implies \mesoreact = \microreact\label{dyn2}\\
\mesodiff > \microdiff &\implies \mesoreact < \microreact\label{dyn3}.
\end{align}
The reaction dynamics is better resolved on the microscopic scale than on the mesoscopic scale, and so we expect the mesoscopic accuracy to increase as $\mesoreact$ approaches $\microreact$ from above. This was shown to be true in \cite{HHP2}. However, as $\mesoreact$ decreases further, the accuracy also worsens. This was also shown in \cite{HHP2}.

Thus, since $\mesodiff$ increases with decreasing $h$, we note that in general we expect the most accurate mesoscopic simulations by selecting $h$ such that \eqref{dyn2} holds. Solving $\mesodiff=\microdiff$ for $h$ yields the optimal mesh size $\hstar$ \cite{HHP2}:
\begin{align}
\hstar = \begin{cases}
\frac{2C_3}{3}\pi\sigma\approx 3.2\sigma, \quad (3D)\\
\sqrt{\pi}e^{\frac{3+2C_2\pi}{4}}\rrad\approx 5.1\sigma, \quad (2D)
\end{cases}
\end{align}
where
\begin{align}
C_3 \approx 1.5164\\
C_2 \approx 0.1951.
\end{align}

\subsection{Hybrid methods}
\label{sec:hybridbackground}

We previously developed a hybrid method \cite{hybrid1} that allowed a given system to be split into two parts: a mesoscopic part and a microscopic part. Species are divided into one subsystem simulated on the microscopic scale, and one subsystem simulated on the mesoscopic scale. The division could depend on spatial constraints, so that a species would be simulated as part of the microscopic subset only in certain parts of space, but not in others. With the system split into two subsets, we would proceed to simulate the system in sequence:
\begin{enumerate}
\item Initialize and set $t=0$. Let the final time be $T$. Select a splitting time step $\Delta t$.
\item Simulate the mesoscopic subset for $\Delta t$ seconds, while keeping the microscopic subset fixed.
\item Simulate the microscopic subset for $\Delta t$ seconds, while keeping the mesoscopic subset fixed. However, we allow microscopic molecules to react bimolecularly with mesoscopic molecules.
\item Synchronize and assign all newly created molecules to their respective scale.
\item Add $\Delta t$ to $t$. Repeat 2-4 until $t=T$.
\end{enumerate}

A crucial and counter-intuitive aspect of this algorithm is that it is necessary to select a time step $\Delta t$ that is neither too small nor too large. The method does not, in general, converge with $\Delta t\to 0$. A rule of thumb is that $\Delta t$ should be selected such that molecules diffuse on the length scale of individual voxels in between synchronization. However, it is straightforward to design a system that would require a smaller $\Delta t$ for accurate simulation, and for which the above hybrid method does not work. 

To see why the simple scheme outlined above leads to the existence of an optimal timestep $\Delta t$, consider the following model system:

\begin{align}
\label{system-main}
\ce{S_1 ->[$k_1$] S_{11} + S_{12} ->[$k_2$] S_2},
\end{align}
with $S_1$ microscopic and $S_{11}$, $S_{12}$, and $S_2$ mesoscopic. When $S_1$ dissociates, $S_{11}$ and $S_{12}$ are placed at contact and they might therefore rebind quickly to form $S_2$. If $\Delta t$ is large, then it is likely that this fast interaction will be captured on the microscopic scale during that time step, and the accuracy will consequently be high. If $\Delta t$ is small on the other hand, then $S_{11}$ and $S_{12}$ will become mesoscopic quickly after $S_1$ dissociates, and information about the spatial correlation of $S_{11}$ and $S_{12}$ will be lost.

Below we propose a way to improve the algorithm to address this problem.

\section{A Convergent Hybrid Method}
\label{sec:method}

Here we propose a hybrid method which builds on the algorithm \cite{hybrid1} outlined in Sect. \ref{sec:hybridbackground} but improves it in two critical ways: First, we propose a new scheme to make the simulation convergent as the splitting time step $\Delta t \to 0$, and second, we use the theory in Sect. \ref{sec:theory} below to enable automatic system partitioning.

\subsection{Algorithm}
\label{sec:hybrid}

To make the method converge monotonically as its time step decreases, we here generalize the splitting over species to allow for \emph{dynamic splitting}, in which a \emph{time-dependent} function maps molecules to either scale. Let $t_j$ denote the time elapsed since the molecule with index $j$ was created, and let $F(S_j,t_j)$ denote the function mapping a molecule of species $S$ of age $t_j$ to either the mesoscopic subset or the microscopic subset. Then, for the system \eqref{system-main} we let
\begin{align}
F(S_1,t) &= \text{microscale}, \quad \text{for all } t\\
F(S_{11},t) &= F(S_{12},t) = \begin{cases}
\text{microscale}, \quad t\leq t_m,\\
\text{mesoscale}, \quad t> t_m,
\end{cases}\label{eq:tm}\\
F(S_2,t) &= \text{mesoscale}, \quad \text{for all } t
\end{align}
where $t_m$ is chosen sufficiently large (see Sect. \ref{selecttm} for how to choose $t_m$), and where $t$ is the time since the molecule was created.

The algorithm in \cite{hybrid1} now becomes a special case ($t_m = 0$) of the algorithm proposed here:
\begin{algorithm}[H]
\caption{\label{newhybrid}\, Hybrid method.}
\begin{enumerate}
\item Initialize the system. Set the time $t=0$. Let $T$ be the length of the simulation.
\item Assign molecules to the mesoscopic and microscopic subsets according to $F(S,t)$.
\item Simulate the mesoscopic molecules for $\Delta t$ seconds. Mesoscopic molecules cannot interact with microscopic molecules during the time step. Any molecules produced will, for the remainder of the time step, be simulated on the mesoscopic scale.
\item Simulate the microscopic molecules for $\Delta t$ seconds, while freezing the mesoscopic molecules. Microscopic molecules can react mesoscopically with mesoscopic molecules. Any molecules produced will, for the remainder of the time step, be simulated on the microscopic scale.
\item Add $\Delta t$ to $t$.
\item Repeat 2-5 until $t=T$.
\end{enumerate}
\end{algorithm}

\subsection{How to split a system}
\label{sec:split_system}

For a given system we will need to determine a suitable splitting in order to achieve high accuracy as well as efficient simulations. Again, consider the system in Eq. \eqref{system-main}. Taking symmetry into account, and by observing that for this particular system the species $S_2$ can be safely simulated on the mesoscopic scale as it only diffuses, we can split the system in five different ways:
\begin{align}
{\cal{X}}_1: \,\, S_1 \text{ micro};\, S_{11}, S_{12}, S_2\, \text{ meso} \\
{\cal{X}}_2: \,\, S_1, S_{11} \text{ micro};\, S_{12}, S_2\, \text{ meso} \\
{\cal{X}}_3: \,\, S_1, S_{11}, S_{12} \text{ micro};\, S_2\, \text{ meso} \\
{\cal{X}}_4: \,\, S_{11}, S_{12} \text{ micro};\, S_{1}, S_2\, \text{ meso} \\
{\cal{X}}_5: \,\, S_{11} \text{ micro};\, S_{1}, S_{12}, S_2\, \text{ meso}
\end{align}
We will now consider the accuracy and convergence of each of the splittings ${\cal{X}}_1-{\cal{X}}_5$.

\underline{${\cal{X}}_1$}: 
In Alg. \ref{newhybrid} the $S_{11}$ and $S_{12}$ molecules will be simulated on the mesoscale for $t_m$ seconds, and therefore, if $t_m$ is chosen large enough such that the molecules either rebind or can be considered well-mixed inside their respective voxels at the end of $t_m$, the system will be accurately simulated. Importantly, this is also true for $\Delta t\to 0$.

\underline{${\cal{X}}_2$}: The accuracy is the same as for the splitting in ${\cal{X}}_1$, since $S_{12}$ is mesoscopic. Association events between microscale and mesoscale molecules have the same spatial resolution as a pure mesoscopic association event.

\underline{${\cal{X}}_3$}: All molecules of interest are simulated on the microscopic scale; the accuracy is therefore the same as for a pure microscale simulation, but with no efficiency gained.

\underline{${\cal{X}}_4$}: The $S_1$ molecule dissociates on the mesoscopic scale, so all spatial correlation is lost (up to the size of the voxel) upon dissociation. Even though we proceed to simulate the $S_{11}$ and $S_{12}$ molecules on the microscopic scale, we still get the accuracy of a mesocopic simulation.

\underline{${\cal{X}}_5$}: The argument from ${\cal{X}}_4$ holds here as well. The accuracy will be the same as for a mesoscopic simulation.

In conclusion, for this model problem the \emph{only} viable splitting of the system \eqref{system-main} (apart from the trivial pure microscopic simulation) is ${\cal{X}}_1$, assuming that we want to simulate as few species as possible on the microscopic scale.

Note that we may, in some cases, be able to simulate the system described by Eq. \eqref{system-main} accurately with the method in \cite{hybrid1}. However, with that algorithm, we cannot take the splitting time step arbitrarily small. The reason for this is described in detail in \cite{hybrid1}. While this is acceptable for some systems, it will lead to inaccuracies for many others.

As an example, consider the following system:

\begin{align}
\ce{S_1 ->[$k_1^1$] S_{11} + S_{12} ->[$k_2^1$] S_2}\label{reac1}\\
\ce{S_2 ->[$k_1^2$] S_{21} + S_{22} ->[$k_2^2$] S_3}\label{reac2},
\end{align}
where $k_2^1$ and $k_2^2$ are large, so that both association reactions are diffusion limited.

First consider the case of $t_m^1 = 0$ and $t_m^2 = 0$. The method now reduces to the method in \cite{hybrid1}. We know that $\Delta t$ has to be large enough.  $S_1$ dissociates on the microscale, so both $S_{11}$ and $S_{12}$ will be microscale until the end of the time step. However, if they are created near the end of the time step, they are likely to survive until the end of the time step, and then turn mesoscopic. If they react on the mesoscale, then the product $S_2$ will be mesoscopic until the end of the time step. If $1/k_1^2\ll\Delta t$ the $S_2$ molecule is likely to dissociate before the end of the time step, in which case $S_{21}$ and $S_{22}$ are initially mesoscopic. Information about the spatial correlation between $S_{21}$ and $S_{22}$ is lost, up to the size of the voxel, and the rest of the simulation may therefore be inaccurate.

With Alg. \ref{newhybrid}, we can guarantee high accuracy, by fixing $t_m^1$ and $t_m^2$ large enough and choosing $\Delta t$ small enough. We then ensure that $S_{11}$ and $S_{12}$ exist on the microscale long enough to either react quickly or become well-mixed insider their respective voxels, while also ensuring that $S_2$ is mesoscopic only on a time scale much shorter than $k_1^2$ (by selecting $\Delta t \ll 1/k_1^2$).

\subsection{Criteria for selecting modeling scale}
\label{sec:theory}

Based on the work outlined in Sect. \ref{mesoparams}, it is clear that choosing a mesh size $h=\hstar$ in general leads to the most accurate mesoscopic simulations. In fact, it is possible to push $h^*$ almost to the size of the molecules if the model is extended to allow for reactions between molecules in adjacent voxels \cite{HP1,FangeSRDME}. The problem is that $\hstar$ is small and this makes the mesoscopic simulations expensive, sometimes significantly more expensive than microscopic simulations \cite{HP1}. Another problem is that $\hstar$ is a function of the reaction radius, and thus different reactions may require different mesh resolutions to be resolved. This can make it impossible to simulate a system accurately with the RDME \cite{HP1}. We therefore want to perform mesoscopic simulations for the majority of the system with $h\gg\hstar$, and handle reactions that require a very fine mesh with a microscopic solver.

For any given $h \ge h^*$,  the RDME solver will match the mean binding time of the microscopic model if the mesoscale propensity functions from \cite{HHP2} are used, but if $h>\hstar$ it will not perfectly match the fine-grained reaction dynamics. Here, we use the relative error in $\mesoreact$ to estimate this error. We let
\begin{align}
\label{eq:W}
W(h) &= \frac{\left| \mesoreact-\microreact \right|}{\microreact}\\
&= \frac{\frac{N}{\krme}-\frac{V}{\krmi}}{\frac{V}{\krmi}}\\
&= \frac{\krmi}{h^d\krme}-1,
\end{align}
where we have used that $\mesoreact>\microreact$ for $h>\hstar$. Now using \eqref{rate-eqs} in place of $\krme$, we obtain
\begin{align}
W &= \frac{\krmi}{h^d\frac{\krmi}{h^d}\left( 1+\frac{\krmi}{D}G(h,\rrad) \right)^{-1}}-1\\
\label{W-expression}
&= \frac{\krmi}{D}G(h,\rrad).
\end{align}

We assume that for $W(h)<\epsilon$, for some small enough $\epsilon$, a reaction is sufficiently resolved on the mesoscopic scale. We can hence use $W(h)$ to decide which species need to be handled on the microscopic scale in order to resolve the reaction dynamics to high enough accuracy. The assumption $W<\epsilon$ with \eqref{rate-eqs} and \eqref{W-expression} holds if and only if
\begin{align}
\frac{\krmi}{h^3}(1+\epsilon)^{-1} < \krme \leq \frac{\krmi}{h^3}.
\end{align}
In other words, a reaction is well resolved mesoscopically when the mesoscopic reaction rate is sufficiently close to the microscopic reaction rate (scaled by the volume of the voxel).

In Sect. \ref{sec:example2} we suggest a reasonable default value for $\epsilon$ based on numerical experiments.

\subsection{How to select $t_m$}
\label{selecttm}

We need to select $t_m$ in Eq. \eqref{eq:tm} so that the hybrid method accurately simulates the system in Eq. \eqref{system-main}. In particular we want the relative error
\begin{align}
\label{hybriderror}
E_{hybrid} = \frac{|\hybridreact-\microreact |}{\microreact}
\end{align}
to be small for $t_m$ large enough.
\noindent
Following a dissociation of $S_1$, the average time until the molecules react is given by
\begin{align}
\label{hr1}
\hybridreact = \microreact|_{t\leq t_m}+\mesoreact|_{t>t_m},
\end{align}
where $\microreact|_{t\leq t_m}$ is the average time until the molecules react on the microscopic scale, given that they react before $t_m$, and $\mesoreact|_{t>t_m}$ is the average time until the molecules react on the mesoscopic scale, given that they react after $t_m$.

Let $S(t)$ denote the probability that the molecules do not react before time $t$. Then
\begin{align}
\label{hr2}
\mesoreact|_{t>t_m} = S(t_m)(\mesodiff|_{t>t_m}+\mesoreact),
\end{align} 
where $\mesodiff|_{t>t_m}$ is the average time that it takes for the molecule to diffuse back to the origin voxel, given that the molecules did not react before time $t_m$. By a similar argument, we can write
\begin{align}
\microreact &= \microreact|_{t\leq t_m} + \microreact|_{t>t_m}\label{mr1} \\
&= \microreact|_{t\leq t_m}+S(t_m)(\microdiff|_{t>t_m}+\microreact).\label{mr2}
\end{align}
Now, by \eqref{hr1}, \eqref{hr2}, \eqref{mr1}, and \eqref{mr2}
\begin{align}
\label{eheq}
E_{hybrid} &= \frac{\hybridreact-\microreact}{\microreact} \\
&= S(t_m)\frac{\mesodiff|_{t>t_m}-\microdiff|_{t>t_m}+\mesoreact-\microreact}{\microreact} \\
&= S(t_m)\left(Q + \frac{\krmi}{D}G(h,\rrad) \right)\\
&\le S(t_m)(Q+\epsilon),
\end{align}
with 
\begin{align}
Q = \frac{\mesodiff|_{t>t_m}-\microdiff|_{t>t_m}}{\microreact}.
\end{align}
It is now easy to see that since we are considering a bounded domain, where $\mesodiff|_{t>t_m}\to \mesodiff$, $\microdiff|_{t>t_m}\to\microdiff$, then as $t_m\to\infty$
\begin{align}
S(t_m)\to 0,\,\, Q \to 0, \,\, \text{as } t_m\to\infty.
\end{align}

The method thus converges as $t_m\to\infty$, which is also easy to see intuitively, as the simulation in practice becomes purely microscopic for $t_m$ large enough.

We have seen numerically that the mesoscopic method incurs an error in average rebind time that is on the order of the time that it takes for a molecule to become well-mixed inside a voxel \cite{HHP2}. We therefore propose to select $t_m=K^2V_{\rm{vox}}^{2/3}/(6D)$, that is, a time proportional to the time that it takes for a molecule to diffuse a distance proportional to the length scale of a voxel. Here $K$ is a constant, and we have found $K=6$ to be sufficiently large.

\begin{figure}
\subfigure{\includegraphics[width=0.49\linewidth]{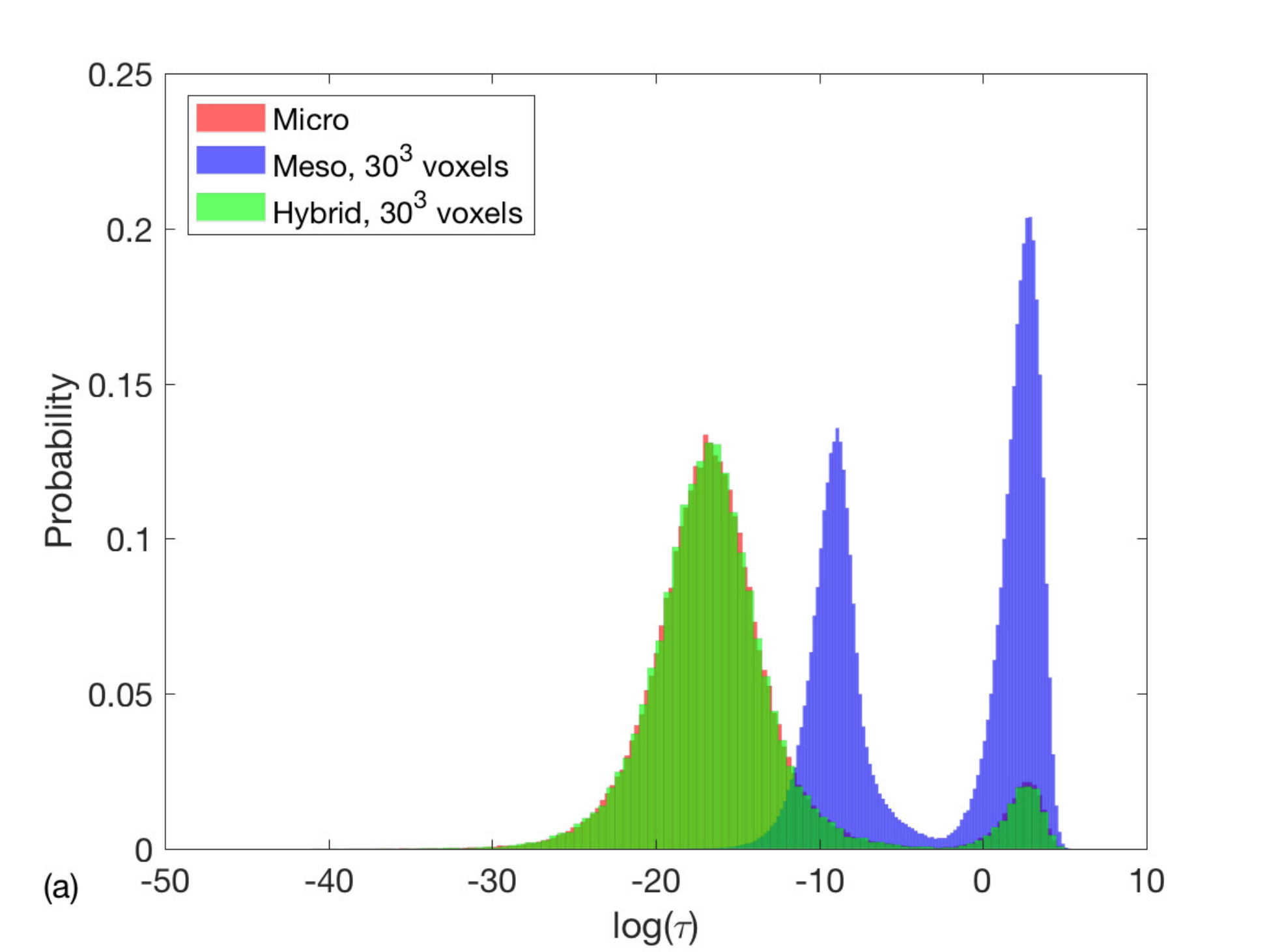}}
\subfigure{\includegraphics[width=0.49\linewidth]{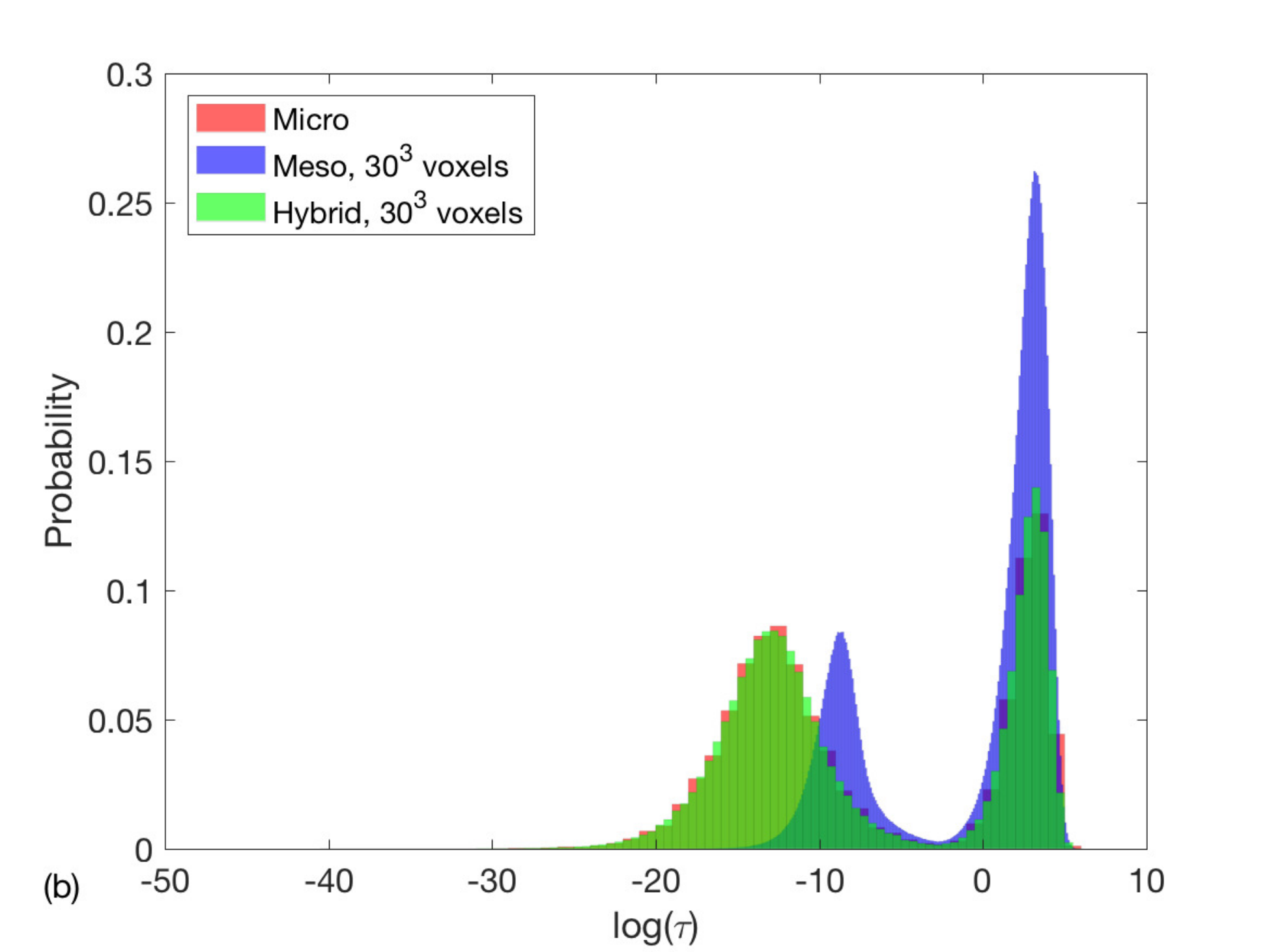}}
\caption{\label{rebindfig} We consider two molecules, one fixed and one diffusing at diffusion rate $1.0$, in a cube of volume $1.0$. The molecules react irreversibly with reaction rate $\krmi$. Let $\tau$ denote the time until the molecules react, given that they start in contact (or in the same voxel on the mesoscopic scale), with a total reaction radius of $0.005$. Above we plot the distribution of the logarithm of the rebind time $\tau$. In (a), $\krmi=1.0$, and in (b), $\krmi=0.1$. We can see that the hybrid method reproduces the microscopic distribution closely (note that green and red overlaps in the figures above), while the mesoscopic RDME does not reproduce the same distribution for diffusion-limited reactions. In particular, the RDME is unable to accurately resolve reaction events occurring on a spatial scale of one voxel or less.}
\end{figure}

We note that $S(t)$ is not known in general since we are considering a bounded domain. However, assuming that $t_m$ will be small enough, and the molecules are some distance from the boundary, we can approximate $S(t)$ by the survival probability for an unbounded domain. This quantity is known analytically, and in 3D is given by \cite{kimshin}:
\begin{align}
\label{survprob}
S(t) = 1-\frac{\krmi}{4\pi\rrad D+\krmi}\left( 1-\exp(\alpha^2t)\rm{erfc}(\alpha\sqrt{t}\right),
\end{align}
where
\begin{align}
\alpha = \left(1+\frac{\krmi}{4\pi\rrad D}\right)\frac{\sqrt{D}}{\rrad}.
\end{align}
We thus know all terms of $E_{hybrid}$ except for $Q$. 
By selecting $t_m$ large enough, we are ensuring that either $Q\approx -\krmi/D G(h,\rrad)$, or that $S(t_m)\approx 0$, and hence $E_{hybrid} \approx 0 $. We show numerically in Fig. \ref{rebindfig} that our choice of $t_m$ is sufficiently large to accurately reproduce the microscopic distribution of rebind times for two different diffusion-limited reaction rates.

\subsection{A more complex example}

To further illuminate how to practically implement the automatic splitting of a system, we will consider some more complex cases. The simple sequence
\begin{align}
\label{simple-system}
\ce{S_1 -> S_{11} + S_{12} -> S_2},
\end{align}
serves as the base case. A single molecule produces two new molecules, spatially correlated, that can then react to form a new molecule. However, we can consider more complex variants of this simple case.

For example, consider the following system:
\begin{align}
\label{splitting-example-1}
\ce{$S_1$ -> $S_{11}$ + $S_{12}$}\\
\ce{$S_{11}$ -> $S_{11}^\ast$}\\
\ce{$S_{11}^\ast$ + $S_{12}$ -> $S_2$}.
\end{align}

In this case $S_1$ produces two molecules, $S_{11}$ and $S_{12}$, that do not directly react. However, if the reaction $\ce{S_{11} -> S_{11}^\ast}$ is fast enough, the above system behaves similarly to the simple system in \eqref{simple-system}. We might therefore need to simulate $S_1$ on the microscopic scale.

We employ the following recursive strategy to identify molecules that should be simulated on the microscopic scale:
\begin{enumerate}
\item Identify all bimolecular reactions for which $W>\epsilon$, for some sufficiently small $\epsilon$, where $W$ is the relative error in the mesoscopic mean binding time, as defined in \eqref{eq:W}.
\item If the two reactants are produced by a dissociating molecule, the dissociating molecule is simulated on the microscopic scale.
\item If they are not, find all reactions producing either of the two molecules, or both.
\item For each sequence of reactions that produces the two reactants, determine whether it starts with a molecule dissociating. If so, that molecule is a candidate to be simulated on the microscopic scale.
\item If not, repeat the process, until we find no new sequences.
\end{enumerate}

For the system \eqref{splitting-example-1} we would therefore first identify that $S_2$ is produced by the molecules $S_{11}^\ast$ and $S_{12}$. These molecules are not produced through any dissociation reaction. We therefore proceed to look for reactions producing one of the two molecules, or both. We find that $S_{11}$ produces $S_{11}^\ast$. We now look for dissociating molecules producing $S_{12}$ and $S_{11}$. We find that $S_1$ produces $S_{12}$ and $S_{11}$, and therefore $S_1$ is a candidate for the microscopic subset.

Another example is the system
\begin{align}
\label{splitting-example2}
\ce{$S_1$ -> $P$ + $S_{12}$}\\
\ce{$P$ -> $D$ + $S_{11}$}\\
\ce{$S_{12}$ -> $S_{12}^\ast$}\\
\ce{$S_{11}$ + $S_{12}^\ast$ -> $S_2$ }
\end{align}

We first find that the only bimolecular reaction is $\ce{$S_{11}$ + $S_{12}$ -> $S_{2}$}$. Assume that we have $W>\epsilon$. We find no dissociation reaction producing $S_{11}$ and $S_{12}$. We now look for reactions producing either of the molecules, or both. No reactions produce both, but we find that $S_{12}$ produces $S_{12}^\ast$, and that $P$ produces $S_{11}$ and $D$. We proceed to look for any reactions producing either $S_{12}$ or $P$, or both. We now find only one such reaction, $\ce{$S_1$ -> $P$ + $S_{12}$}$. Since this is a dissociation reaction, $S_1$ will be simulated on the microscopic scale.

\section{Implementation}
\label{sec:implementation}
We have implemented the hybrid method as an extension to the high-level PyURDME \cite{molns} Python API. This allows for specification of microscopic/hybrid systems and execution via a simple, object-oriented Python modeling interface. In the following sections, we describe the different components of the solver and discuss computational complexity and performance aspects of hybrid simulation.

\subsection{Mesoscopic solver: Next Particle Method}
In the NSM, reaction and diffusion events in each voxel are grouped, and the heap is organized so that each leaf corresponds to a voxel. In each iteration, the next reaction or diffusion event is executed and the next event time is updated along with the heap for each affected voxel. For a fine mesh, the vast majority of events are diffusion events and the simulation cost is dominated by the time to execute diffusion events. Ignoring reactions, the simulation cost, $C_{NSM}$, on a uniform grid with $N$ voxels and $M$ molecules of a single diffusing species can be written as
\begin{align}
C_{NSM}(N,M) = C_1 N^{2} M\log{N},
\end{align}
where $C_1$ is an implementation- and architecture-dependent constant. Here we instead propose a particle-centric mesoscopic algorithm, the Next-Particle Method (NPM), that tracks individual particles on the grid. We simulate a mesoscopic system as follows:
\begin{enumerate}
\item Particles are stored in a list, with information about species type and which voxel they currently occupy.
\item For each particle we generate the time to and the destination voxel of its next diffusion event. Add each diffusion event to the heap. For $N$ particles, the size of the heap will be $N$.
\item For two reactive particles occupying the same voxel, we generate the time until the next tentative event, and add that event to the heap.
\item Execute the next event.
\item Update all dependent events. If a molecule diffused, add its next diffusion- and reaction events to the heap. If molecules reacted, add new diffusion and reaction events to the heap for all molecules that were affected.
\item Repeat 4-5 until the end of the simulation.
\end{enumerate} 
The main advantage of the NPM in the context of the hybrid method is that it minimizes the overhead of switching between the mesoscopic and the microscopic solvers, since the two solvers can share one datastructure for the particle list. The microscopic solver needs to know the position of individual molecules, so maintaining one particle list simplifies the mapping between discrete positions on a grid and continuous positions in space.
The cost for the NPM for the example above can be written
\begin{align}
C_{NPM}(N,M) = C_2 N^{2} M\log{M},
\end{align}
\noindent
where $C_2$ is a constant. This highlights another potential advantage of the NPM in the context of hybrid simulation: it can be more efficient than NSM for highly resolved meshes if the number of voxels are larger than the number of particles. This will often be the case for highly resolved geometries. \

Note however, that the hybrid framework proposed here does not depend on the particular implementation of the mesoscopic method, and that it would be possible to alternate algorithms depending on the particular values of $N$ and $M$.

\subsection{Microscopic solver: GFRD}
\label{micro-implementation}
On the microscopic scale, the system is simulated with the Smoluchowski model, as described in Sect. \ref{micro-background}. For a system of more than one or two molecules, we have an intractable many-body problem. To deal with this, we employ a strategy conceptually similar to the GFRD algorithm \cite{ZoWo5a}.

The first step of the algorithm is to divide the system into subsets of one or two molecules, and to select a time step $\Delta t$, such that to high accuracy we can update the subsets independently during $\Delta t$. Molecules that are each others' nearest neighbors are updated in pairs, while all other molecules are updated as single molecules. The time step $\Delta t$ is chosen as large as possible, with the constraint that the probability of interactions between the separate subsets is small. We then simulate each subset for $\Delta t$ seconds. 

For each subset we look for the next reaction: if two molecules react bimolecularly, we can sample the time until they react from $p_{\fatr}(\ast,t|\fatr_0,t_0)$ (defined in \eqref{micro_time}), if either or both molecules can react unimolecularly, we can sample the next reaction time from exponential distributions, and finally we will look for possible interactions with the boundary. The reaction that occurs first is executed. We then repeat the procedure until the subsystem has been advanced to time $t_0+\Delta t$.

Instead of solving Eq. \eqref{smolu-eq1} with boundary conditions given by Eqs. \eqref{smolu-eq1-bicond1}, \eqref{smolu-eq1-bicond2}, and \eqref{smolu-eq1-bicond3} exactly, we solve it using the operator split approach described in \cite{SHeLo11}. Furthermore, we only sample from $p_{\fatr}(\ast,t|\fatr_0,t_0)$ and $p_{\fatr}(\fatr,t|\fatr_0,t_0)$ if the distance between the molecules is small and the probability of a reaction is fairly large. If the probability of reaction during the time step is small, it can be more efficient to simulate the pair using brute-force Brownian Dynamics until the molecules are close. The cut-off is typically at a distance of around a few reaction radii.

\subsection{Hybrid solver}
\label{hybrid-complexity}
We can now couple the mesoscopic NPM with the microscopic solver in a simple loop. Since the NPM keeps track of individual molecules, it is straightforward to map each molecule to either scale according to the splitting function $F$. Both solvers can keep track of how long a molecule has existed, thus making it easy to determine whether a microscopic molecule can be mapped to the mesoscopic scale.

When a molecule switches from the mesoscopic scale to the microscopic scale, we need to know its position in continuous space. We sample its position from a uniform distribution on its voxel. Similarly, when a molecule switches from the microscopic scale to the mesoscopic scale we need to know which voxel the molecule occupies. This is straightforward, since we track which voxel a molecule occupies to accurately simulate its interaction with the boundary. This process is described in detail in \cite{hybrid1}. The overhead from this switching is inversely proportional to the splitting timestep, $C_{coupling} = C_3 (\Delta t_s)^{-1}$, where we have assumed that the number of particles that switch in each timestep is small, compared to the total number of particles on both scales. 

With $M_1$ the average number of mesoscopic particles and $M_2$ the average number of microscopic particles during the course of a simulation, the complexity of the overall hybrid method can be described by:

\begin{align}
C_{hybrid}(N,M_1,M_2) = \frac{C_2 M_1\log{M_1}}{N^{-2}} + C_{gfrd}(M_2) + C_3 (\Delta t_s)^{-1} .
\end{align}
\noindent
Since the number of particles handled on the different scales depends on the mesh resolution, i.e. $M_1$ and $M_2$ are functions of $N$, the cost of the solver is complicated to estimate a priori, and it implies the existence of an optimal choice of $N$ for performance. This will be illustrated in Sect.~\ref{seceff}.

Note that each solver could be optimized depending on the system. If we were mainly interested in simulating systems in a cube, the microscopic solver could be significantly optimized by simplifying the process of keeping track of the boundary. For a system with many more particles than voxels, we could choose to simulate the mesoscopic part of the system with the NSM rather than with the NPM.

In Sect. \ref{seceff} we show how the contribution to the total execution time of each solver depends on the mesh size and the system. The total cost of a simulation depends non-linearly on the mesh size, since we need to balance the trade-off between a coarse mesh and fast mesoscopic simulations but expensive microscopic simulations, with a fine mesh on which the mesoscopic simulations are more expensive while the microscopic simulations will be faster (due to the fact that we will simulate fewer molecules on the microscopic scale on a fine mesh).

\section{Numerical examples}
\label{sec:results}

While the theory above is derived under the assumption of a Cartesian mesh, we have shown that in most cases it can be applied also to the case of unstructured meshes \cite{rdmeunstruc}, by substituting the voxel width $h$ for $V_{vox}^{1/3}$, where $V_{vox}$ is the volume of a voxel in a mesh. In particular, we show in Sect. \ref{sec:example2} below that we can accurately split and simulate a system on an unstructured mesh.

Furthermore, we demonstrate that we can accurately simulate a problem previously shown to be intractable with the standard RDME model \cite{rdmeunstruc}, and finally we show the existence of an optimal mesh size, from an efficiency perspective, in between the coarsest and finest possible mesh sizes.

\subsection{Splitting Species: Accuracy and Efficiency}
\label{sec:example2}

In this example we demonstrate that for a given system, we can split the species into a microscopic subset and a mesoscopic subset using $W(h)$ defined in Eq. \eqref{eq:W}. The resulting splitting of species should yield accurate and efficient simulations on a given unstructured mesh.

First we want to determine a suitable $\epsilon$ such that $W<\epsilon$ indicates that the reaction is sufficiently resolved on the mesoscopic scale. We again consider the simple system
\begin{align}
\ce{S_1 ->[$k_1$] S_{11} + S_{12} ->[$k_2$] S_2}.
\label{example-simple-system}
\end{align}

In Sect. \ref{sec:theory} we found that
\begin{align}
W = \frac{k_2}{D}G\left(V_{\rm{vox}}^{\frac{1}{3}},\rrad\right)
\end{align}
is a measure of how well the rebinding time of a pair of molecules is resolved (where $k_2^{\rm{meso}}$ is given by \eqref{rate-eqs}). 

While we have no theory relating $W$ directly to the error in the mesoscopic simulation of the system, we can use it as an indirect measure of the error. Via numerical simulations we can find an $\epsilon$ such that $W<\epsilon$ implies that the simulations will be accurate.

For simplicity we consider the system \eqref{example-simple-system} in a cube. We let the microscopic parameters be given by
\begin{align}
\begin{cases}
\sigma_1 = \sigma_{11} = \sigma_{12} = \sigma_2 = 0.0025\\
D_1 = D_{11} = D_{12} = D_2 = 1.0\\
k_1 = 10.0\\
V = (50\hstar\sigma)^3,
\end{cases}
\end{align}
and we sample $k_2$ from $[0.001,1.0]$.
By design we expect the best agreement between mesoscale and microscale simulations for a mesh of $50^3$ voxels. Note that these parameters are chosen arbitrarily, but we will proceed to show that the results can be applied successfully to a numerical example with different parameters.

First we compare pure mesoscopic and microscopic simulations. Let $\mathbf{y_{\rm{me}}} = (y_{\rm{me}}^1,\ldots,y_{\rm{me}}^N)$ be the average number of $S_2$ molecules, computed from $M_{\rm{me}}$ mesoscale trajectories sampled at the time points $t_1,\ldots,t_L$, and let $\mathbf{y_{\rm{mi}}} = (y_{\rm{mi}}^1,\ldots,y_{\rm{mi}}^L)$ be the average number of $S_2$ molecules computed as the average of $M_{\rm{mi}}$ microscale trajectories sampled at the time points $t_1,\ldots,t_L$. We consider the max-norm error E, defined as
\begin{align}
\label{Edef}
E = \max_{1\leq i \leq L}\left| y_{\rm{me}}^i-y_{\rm{mi}}^i \right|.
\end{align}

\begin{figure}
\subfigure{\includegraphics[width=0.45\linewidth]{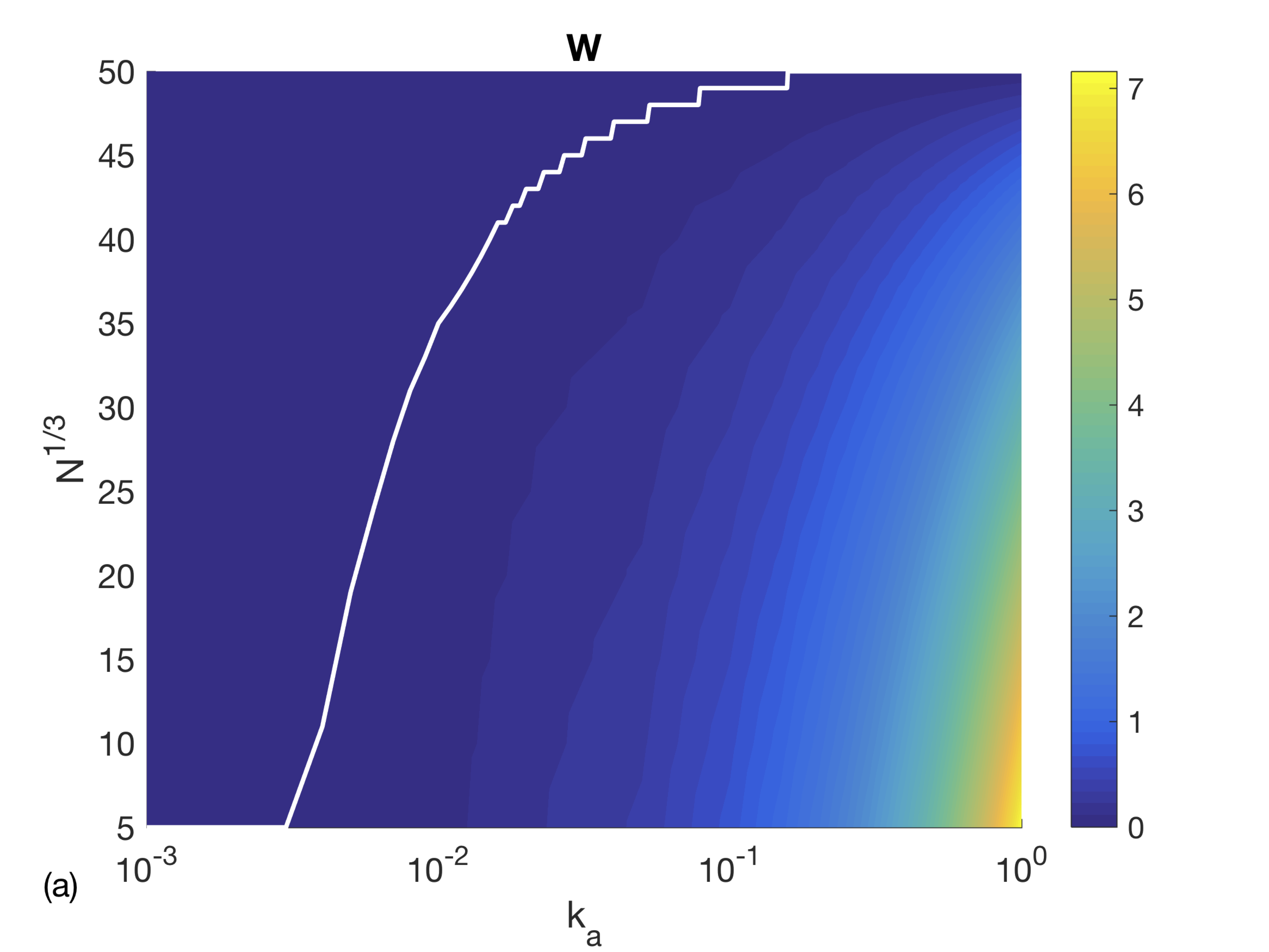}}
\subfigure{\includegraphics[width=0.45\linewidth]{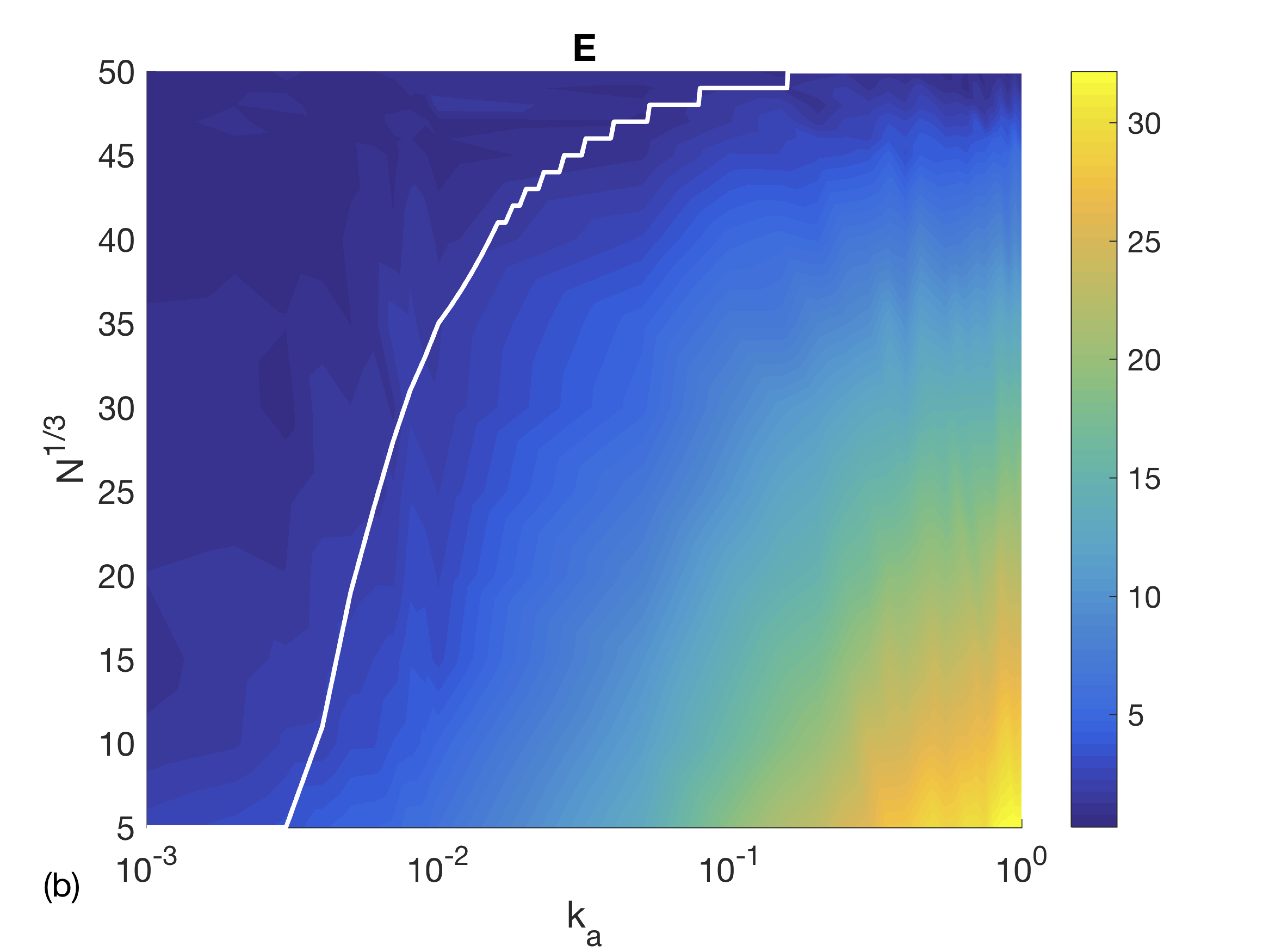}}
\caption{\label{example2_fig1}In (a) we plot $W$ as a function of $N^{\frac{1}{3}}$ and $\krmi$. For points above the white solid line we have $W<0.025$. In (b) we plot the error $E$ as a function of $N^{\frac{1}{3}}$ and $\krmi$. Again, for points above the white solid line, $W<0.025$. For small $\krmi$, the error is small for all mesh sizes, while as $\krmi$ increases, we need a large $N$ to keep the error small. We can see that for $W<0.025$, the error is roughly on the order of 1.}
\end{figure}

For small values of $k_2$ we expect the mesoscopic simulations to be accurate also for coarse meshes, while for large $k_2$, we expect the error to be large unless the spatial resolution is near the maximum resolution of $50^3$ voxels. In Fig. \ref{example2_fig1} we show how $W$ correlates with the error $E$ for different $k_2$, and that $\epsilon=0.025$, although arbitrary, is a reasonable choice yielding an error of the order of 1.

We now apply the choice of $\epsilon=0.025$ to an expanded system of three bimolecular reactions:
\begin{align}
\label{example2-full-system}
\ce{S_1 ->[$k^1_1$] S_{11} + S_{12} ->[$k^1_2$] S_2}\\
\ce{S_2 ->[$k^2_1$] S_{21} + S_{22} ->[$k^2_2$] S_3}\\
\ce{S_3 ->[$k^3_1$] S_{31} + S_{32} ->[$k^3_2$] S_4},
\end{align}
with parameters different from the simple system above. The system is simulated inside a sphere of radius $0.5$, discretized with an unstructured mesh consisting of 6395 voxels.

Depending on the values of $k^i_2$, $i=1,2,3$, we will simulate some combination of $S_1$, $S_2$, and $S_3$ on the microscopic scale. For $W_i>\epsilon$, $S_i$ is simulated on the microscopic scale. The minimum time that a molecule has to exist on the microscopic scale before it becomes mesoscopic is given by $t_m = \frac{V_{\rm{vox}}^{2/3}}{C^26D_i}$, with $C=6$, cf. Sect. \ref{selecttm}. 

We consider six different combinations of reaction rates, see Table \ref{ex1-params-tab}. In each case we will have a different combination of $S_1$, $S_2$, and $S_3$ on the microscopic scale. For $k_2^\ast>0.1$, $W>\epsilon$, while for $k_2^\ast = 0.001$ we have $W<\epsilon$. Thus, for case 6 the hybrid method will simulate all molecules on the mesoscopic scale, and we therefore expect the mesoscopic simulation to agree well with the microscopic simulation.
\begin{table}
\begin{tabular}{ r | c | c | c | c | c | c |}
 & Case 1 & Case 2 & Case 3 & Case 4 & Case 5 & Case 6\\ \cline{2-7}
$k_2^1$ & 0.1 & 0.001 & 0.1 & 0.1 & 0.001 & 0.001\\ \cline{2-7}
$k_2^2$ & 0.1 & 0.3 & 0.3 & 0.001 & 0.001 & 0.001\\ \cline{2-7}
$k_2^3$ & 0.1 & 0.001 & 0.001 & 0.2 & 0.2 & 0.001\\ \cline{2-7}
\end{tabular}
\caption{\label{ex1-params-tab}Association rates for the six different cases.}
\end{table}

\begin{figure}
\subfigure{\includegraphics[width=0.49\linewidth]{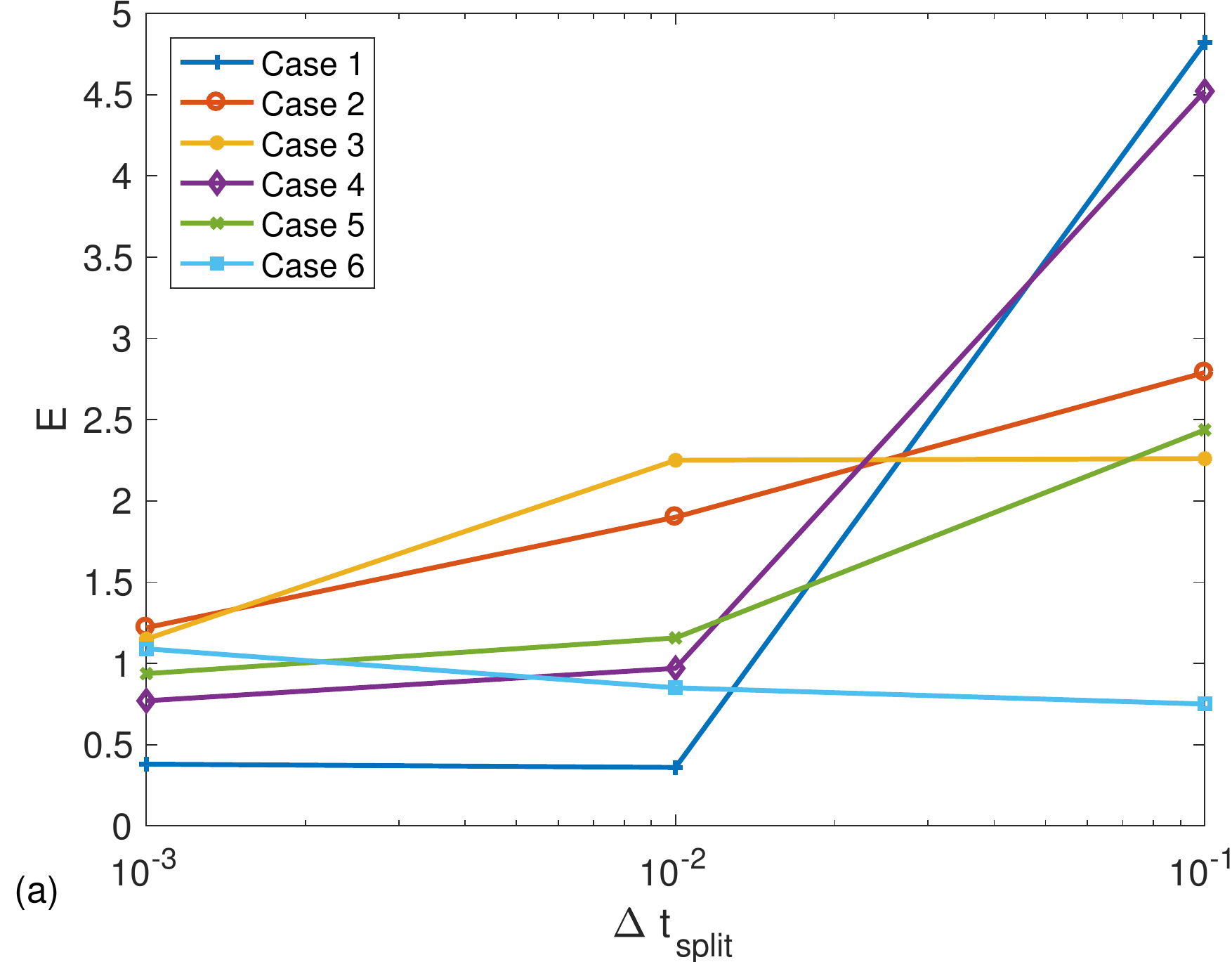}}
\subfigure{\includegraphics[width=0.49\linewidth]{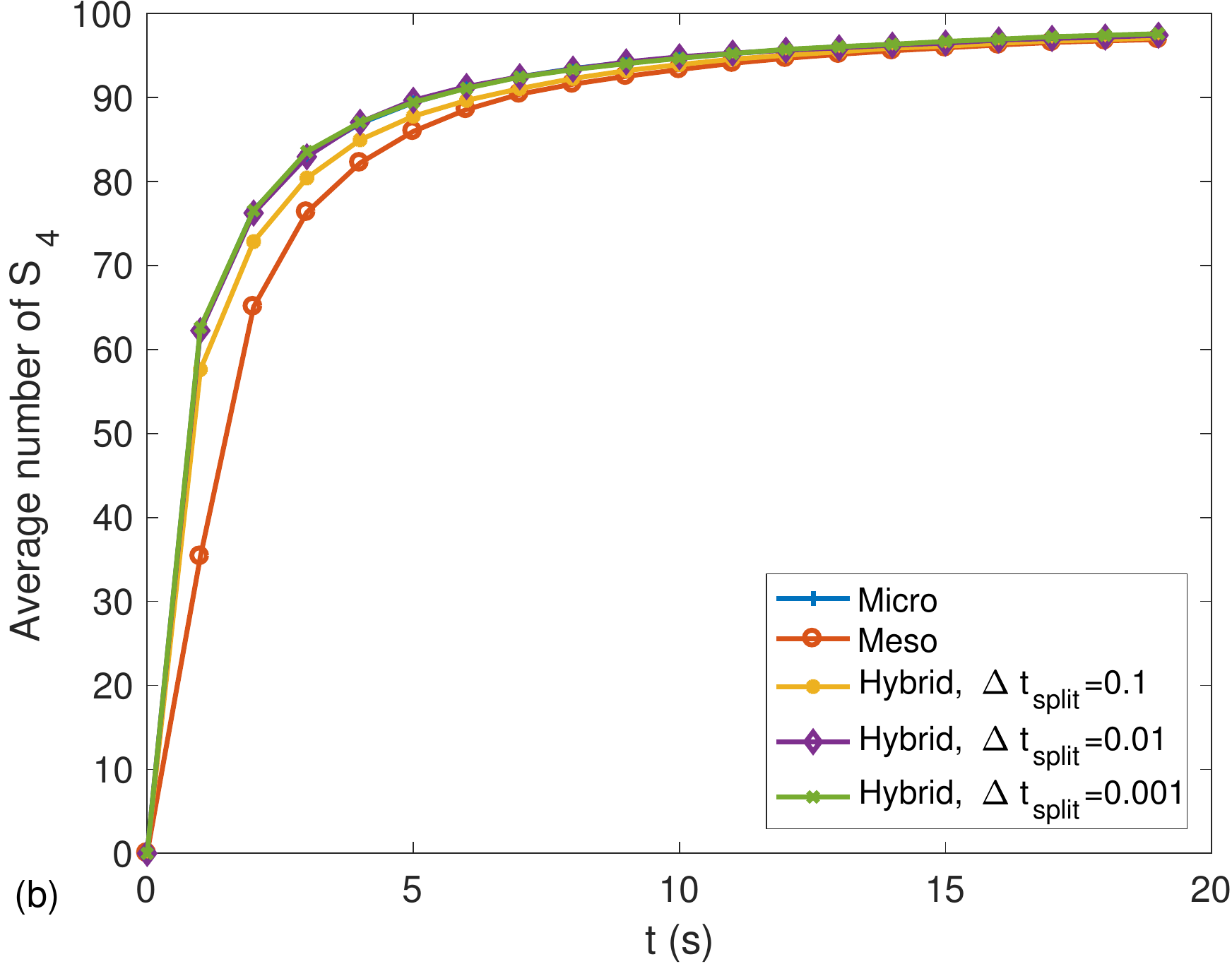}}
\caption{\label{example2_fig2}Left (a): We plot the error $E$, as defined in Eq. \eqref{Edef} as a function of the splitting time step $\Delta t_{\rm{split}}$. For the smallest time step, $\Delta t_{\rm{split}}=0.001$, the error remains small for Case 6, but is larger for Cases 1 and 4, in which some or all of the reactions are diffusion limited. Right (b): The average number of $S_4$ molecules over time in Case 1. We see that the hybrid method with $\Delta t_{\rm{split}}=0.1$ underestimates the average number of $S_4$ molecules, but still produces better results than with a pure mesoscopic simulation. The hybrid method matches the microscopic results closely for $\Delta t_{\rm{split}}\leq 0.01$.}
\end{figure}

In Fig.~\ref{example2_fig2} we show that the hybrid method agrees well with the microscopic simulations. In case 6, the RDME agrees well with the microscopic model, as expected. In addition, we show that the accuracy increases with a decreasing splitting time step. For a relatively large splitting time step of $0.1$, the method produces results with a fairly large error, but as we refine the splitting time step, the results approach that of a pure microscopic simulation. We have tabulated the errors of the hybrid method errors in Table \ref{tab1-error} with the errors of pure mesoscopic simulations. Even with a fairly large splitting time step, the error in the hybrid method is smaller. For case 6, in which the reactions are slow compared to diffusion, both the hybrid method and the RDME produce accurate results.

While one reason to use a hybrid method is to gain efficiency over a very fine-grained RDME simulation, another is that some systems cannot be simulated accurately with a standard RDME model. In \cite{rdmeunstruc} we considered the following system:
\begin{align}
\label{systemnonlocal}
\ce{S_1 ->[$k_d$] S_{11} + S_{12} ->[$k_r$] S_2}\\
\ce{S_2 ->[$k_d$] S_{21} + S_{22} ->[$k_r$] S_3},
\end{align}
where $k_d = 10.0$ and $k_r = 0.1$. If $\sigma_i$ is the reaction radius of species $S_i$, and $\sigma_{ij}$ the reaction radius of species $S_{ij}$, then $\sigma_1 = 10^{-3}$, $\sigma_{11} = 0.8\times 10^{-3}$, $\sigma_{12} = 0.8\times 10^{-3}$, $\sigma_{2} = 2.0\times 10^{-3}$, $\sigma_{21} = 1.8\times 10^{-3}$, $\sigma_{22} = 1.8\times 10^{-3}$, and $\sigma_{3} = 2.5\times 10^{-3}$. For simplicity, all molecules diffuse with diffusion rate $1.0$. The domain is a cube of volume $1.0$. 

To resolve the first association we need a mesh size of around $h^*_1 = \frac{2}{3}C_3\pi(\sigma_{11}+\sigma_{12})\approx 5.0\times 10^{-3}$, and to resolve the second association we need a mesh size of around $h^*_1 = \frac{2}{3}C_3\pi(\sigma_{21}+\sigma_{22})\approx 1.14\times 10^{-2}$. We showed in \cite{rdmeunstruc} that we cannot resolve both reactions simultaneously with the standard local RDME; we could simulate the system by allowing reactions between neighboring voxels. However, these simulations become expensive as the mesh needs to be highly refined, and they cannot be trivially extended to unstructured meshes.

\begin{figure}
\subfigure{\includegraphics[width=0.49\linewidth]{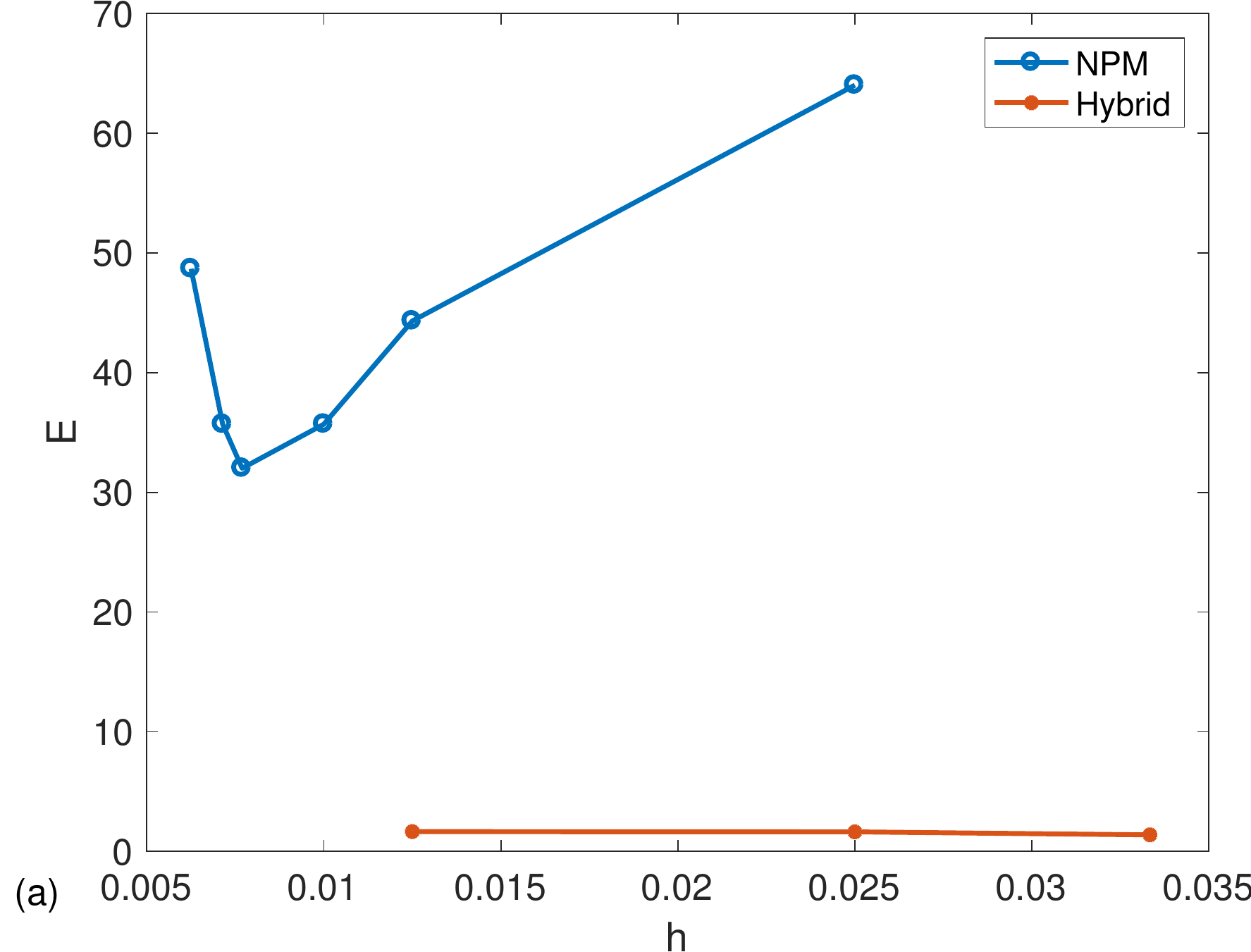}}
\subfigure{\includegraphics[width=0.49\linewidth]{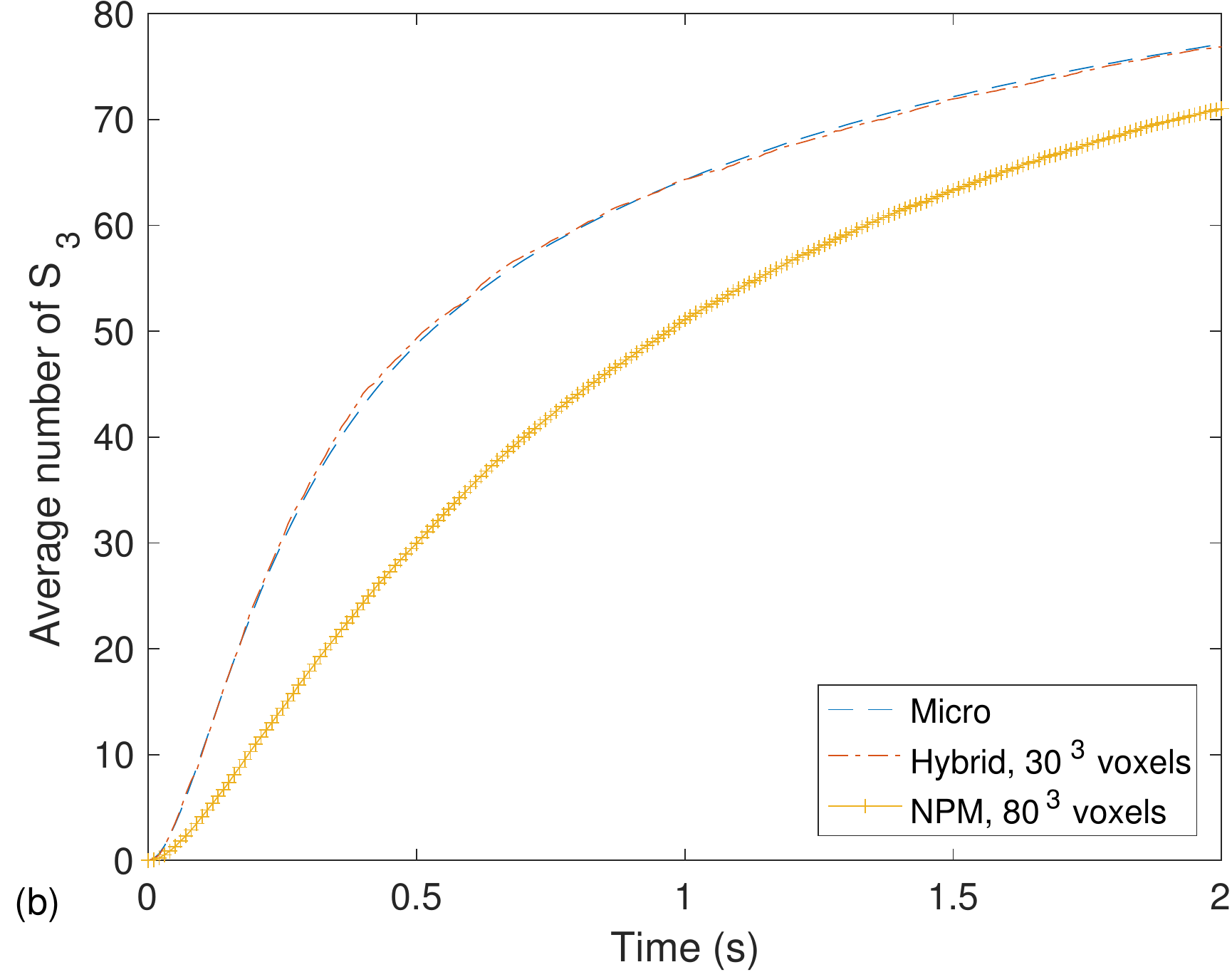}}
\caption{\label{example1nonlocal}Left (a): The error $E$, as defined by Eq. \ref{meanerror}. The RDME does not match the microscopic dynamics for any mesh size. The hybrid method is able to capture the dynamics of the system by simulating the $S_1$ and $S_2$ species on the microscopic scale, and all other species on the mesoscopic scale. Right (b): Example trajectory of the average population of $S_3$ molecules. The hybrid method agrees well with the microscopic results, while the NPM on a mesh of $80^3$ voxels does not agree with the microscopic simulations.}
\end{figure}

We show here that another viable approach is to simulate the system with a hybrid method. The system is simulated for 2 seconds, with 201 uniform time samples including $t=0$. In Fig. \ref{example1nonlocal} we plot the error $E$ as a function of the mesh size, where $E$ is defined as in \cite{rdmeunstruc}. Let ${\cal{S}}=\{S_1,S_{11},S_{12},S_{2},S_{21},S_{22},S_3\}$. Then
\begin{align}
\label{meanerror}
E(h) = \frac{1}{201}\sum_{i=1}^{201} \sum_{S\in{\cal{S}}}| S^{\ast}_{h,i}-S_i^{micro}|,
\end{align}
where $S_i^{micro}$ is the average population of species $S$ at time $t_i$, obtained with the microscopic algorithm, and where $S_{h,i}^{\ast}$ is the average population of species $S$ at time $t_i$ obtained with either the hybrid algorithm or with the NPM, simulated on a mesh with a voxel width of $h$.

\begin{table}
\begin{tabular}{ r | c | c | c | c | c | c |}
 & Case 1 & Case 2 & Case 3 & Case 4 & Case 5 & Case 6\\ \cline{2-7}
RDME & 27.08 & 3.78 & 7.19 & 9.87 & 4.597 & 0.81\\ \cline{2-7}
Hybrid, $\Delta t_{\rm{split}}=0.1$ & 4.82 & 2.79 & 2.26 & 4.52 & 2.437 & 0.75\\ \cline{2-7}
Hybrid, $\Delta t_{\rm{split}}=0.01$ & 0.36 & 1.90 & 2.25 & 0.97 & 1.157 & 0.85 \\ \cline{2-7}
Hybrid, $\Delta t_{\rm{split}}=0.001$ & 0.38 & 1.22 & 1.15 & 0.77 & 0.937 & 1.09 \\ \cline{2-7}
\end{tabular}
\caption{\label{tab1-error} Max-norm error. We see that the hybrid method, even for a large splitting time step, produces results that are more accurate than pure mesoscopic simulations. In the case where all reactions are slow compared to diffusion, the RDME produces accurate results, as does the hybrid method.}
\end{table}

\subsection{Efficiency: Non-linear dependence on the mesh size}
\label{seceff}
As already discussed in Sect. \ref{hybrid-complexity}, the total execution time is the sum of the time spent on the mesoscopic scale, $T_{meso}$, the microscopic scale, $T_{micro}$, and the overhead incurred from the coupling of the scales. The time spent on the microscopic scale depends on how many of the species are microscopic, which in turn depends on the resolution of the mesh. On a fine mesh, we will simulate fewer molecules on the microscopic scale, and for a shorter time, but we pay the price of a more costly mesoscopic simulation.

In this numerical example we show that the total execution time $T$ is a non-linear function of $T_{micro}$ and $T_{meso}$, and that to optimize $T$ we need to balance $T_{micro}$ and $T_{meso}$ in a non-trivial way.

\begin{figure}
\subfigure{\includegraphics[width=0.49\linewidth]{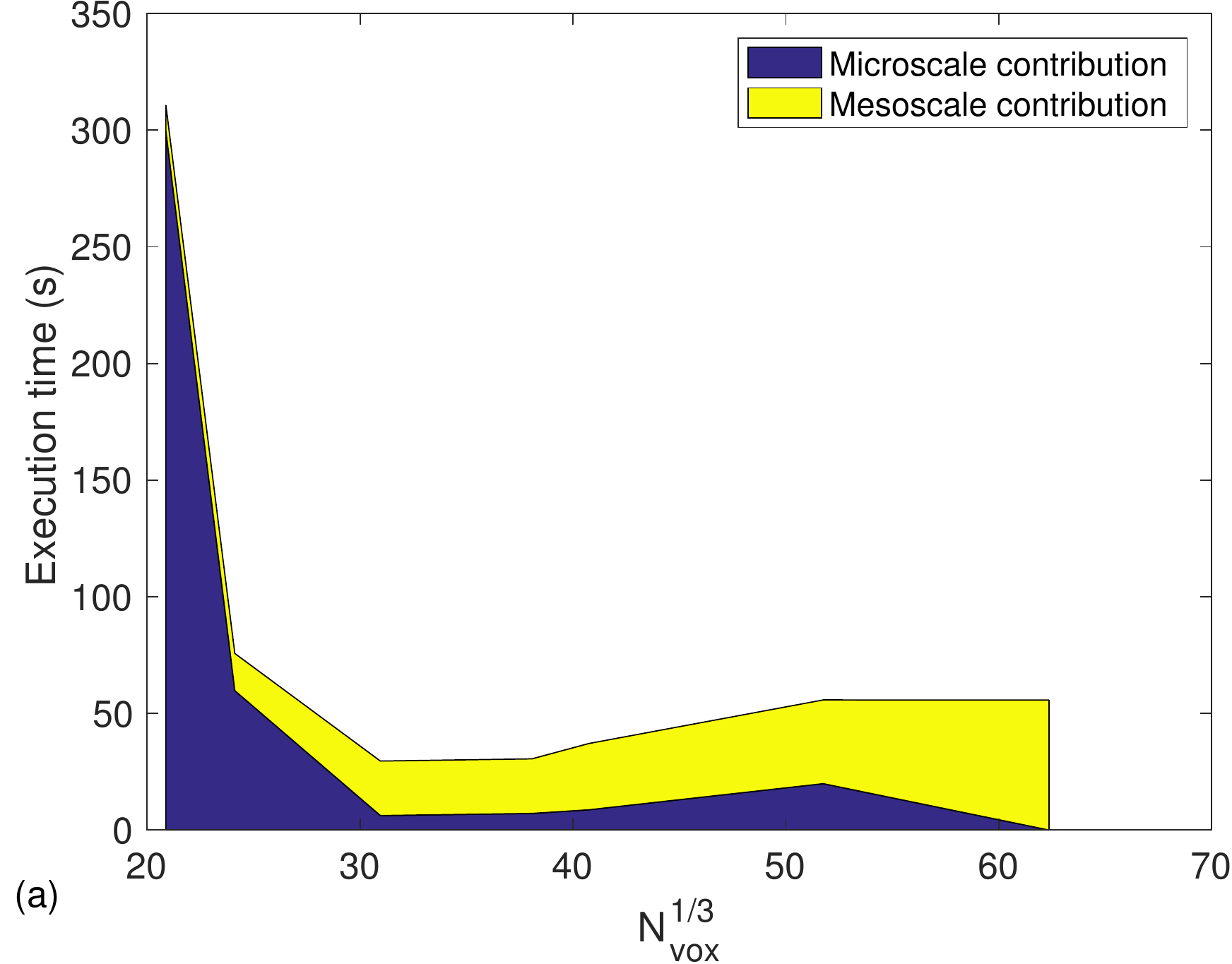}}
\subfigure{\includegraphics[width=0.49\linewidth]{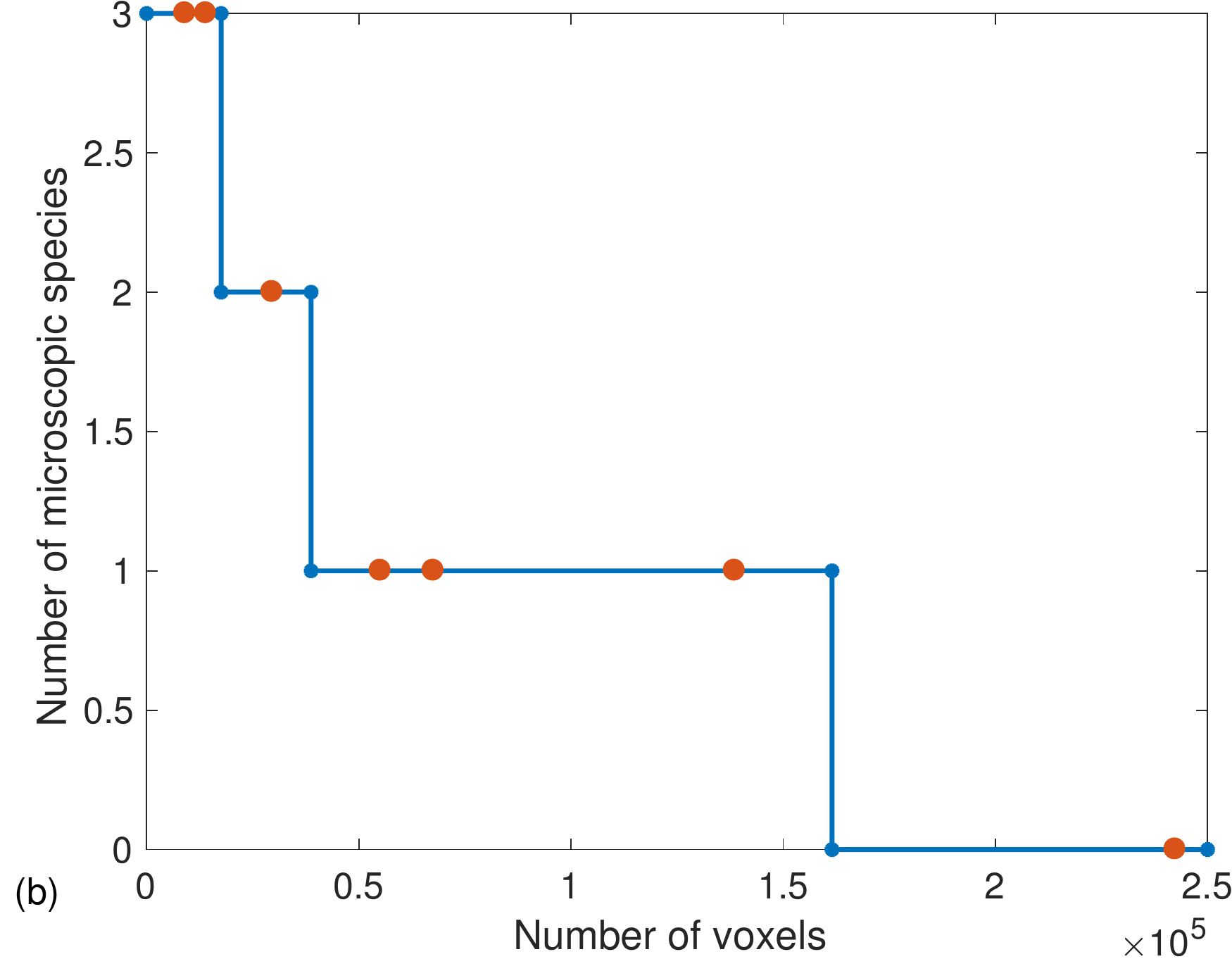}}
\caption{\label{meshsweep-fig}Left (a): The total execution time is a nonlinear function of the mesh size. For a coarse mesh, most of the simulation time will be spent on the microscopic scale. As the mesh is successively refined, the time spent on the mesoscopic scale starts to dominate. The shortest total simulation time is obtained for a mesh of around 30,000 voxels. Right (b): Number of particles on the microscopic scale as a function of the mesh resolution. Red dots indicate the sample points used to generate the plot in (a). As we can see, for the two left-most points, we simulate $S_1$, $S_2$ and $S_3$ on the microscopic scale. As we move to the right, we will simulate two, one, and finally, for the right-most point, no species on the microscopic scale. We can see in (a) that a pure mesoscopic simulation is slower than a simulation with one microscopic particle, but with a much coarser mesh.}
\end{figure}

We consider the system
\begin{align*}
\label{meshsweep-system}
\ce{S_1 ->[$k^1_1$] S_{11} + S_{12} ->[$k^1_2$] S_2}\\
\ce{S_2 ->[$k^2_1$] S_{21} + S_{22} ->[$k^2_2$] S_3}\\
\ce{S_3 ->[$k^3_1$] S_{31} + S_{32} ->[$k^3_2$] S_4}
\end{align*} 
where $k^1_i = 20.0$, $i=1,\ldots,5$ and $k^1_2=0.0016$, $k^2_2=0.00145$, and $k^3_2=0.0014$. We initialize the system with 200 $S_1$ molecules and 200 $S_2$ molecules, all with uniformly sampled positions on the domain. The domain is a sphere with radius 0.5, and we consider a sequence of meshes, ranging from coarse to fine. In Fig. \ref{meshsweep-fig} we show that there exists an optimum, with respect to total execution time, between the coarsest and the most resolved mesh. Note that this particular system can be accurately simulated on any mesh resolution with the hybrid method, and therefore the error remains small for all mesh sizes, with the only thing changing being the number of molecules simulated on either scale, and for how long the microscopic molecules remain microscopic.

In general there is no way to \emph{a priori} determine the optimal mesh size, as it will depend on the system under study as well as the initial condition. It will also depend on the size of the simulation; if we are planning on running many or very long trajectories, making the total simulation time substantial, we can afford an expensive preprocessing step. On the other hand, if the total simulation time is short to moderate, an expensive preprocessing step will not be worthwhile. We therefore propose a heuristic approach to selecting the mesh size.

If we can afford an expensive preprocessing step, we can simulate either full trajectories or shortened trajectories on a sequence of mesh resolutions to find a mesh resolution that appears to minimize the total execution time (there is of course no guarantee that we have found an actual minima). To speed this process up we could, if the system allows it, perform the simulations on a structured Cartesian grid on a regular domain. That way we avoid the costly process of generating a sequence of unstructured meshes. While there is no guarantee that the system behaves the same way on a structured Cartesian mesh as on the actual domain of interest, we can still expect to get an approximation of the relative cost of simulations on different mesh sizes.

Also note that in many cases the mesh size will be constrained by the geometry of the problem. Internal structures could require a certain minimum mesh resolution, meaning that we cannot select the mesh size that optimizes the execution time, but that we instead are constrained to a certain mesh, and have to choose the best splitting given the mesh resolution.

\section{Discussion}
\label{sec:conclusions}

We have developed a hybrid method coupling simulation of the mesoscopic and microscopic modeling scales. The method can, for a certain class of systems, automatically propose a splitting of species based on how diffusion-limited the reactions are. Furthermore, we show that the new method converges with decreasing splitting time step for a larger class of systems than a previously developed method \cite{hybrid1}.

We apply the method to a numerical example, showing that it accurately, and with increased efficiency compared to microscopic simulations, splits the system. We also show how the optimal splitting can be found for a mesh between the coarsest and the finest possible resolutions. It is therefore necessary to find a balance between how many molecules to simulate on the microscopic scale, and how fine the mesh should be.

The approach described in this paper can, in general, be applied to systems where molecules are created in spatial proximity through some sequence of unimolecular and bimolecular reactions.
Another possibility is that molecules are created in spatial proximity due to more complex interactions with internal membranes or fibers; processes not necessarily captured by the scheme outlined above. It is also plausible that microscale resolution could be needed for other reasons, such as for processes where molecules in 3D react with complex membranes or move due to active transport. Automatic splitting of such systems would require a different analysis.

\section{Acknowledgement}
This work has been funded by National
Institute of Health (NIH) NIBIB Award No. R01-EB014877, Department of Energy (DOE) Award No. DE-SC0008975, the Swedish Research Council under Award No. 2015-03964 and the eSSENCE strategic collaboration on eScience.

\end{document}